\documentclass[12pt]{amsart}

\usepackage[frak=euler]{mathalfa}
\usepackage[utf8]{inputenc}
\usepackage[a4paper,top=0.5in,bottom=0.2in,left=1.3in,right=1.3in,footskip=0.3in,includefoot]{geometry}
\usepackage[english]{babel}
\usepackage{amsmath,amssymb,amsfonts,amsthm}
\usepackage{mathrsfs,mathtools,nicefrac}
\usepackage{graphicx}
\usepackage[dvipsnames]{xcolor}
\usepackage{framed,accents,enumerate,comment,wrapfig}
\usepackage{esint}
\IfFileExists{tensor.sty}{\usepackage{tensor}}{}
\usepackage{tikz,tikz-cd}
\usetikzlibrary{calc,angles,quotes}
\usepackage{aliascnt}
\usepackage{hyperref}
\hypersetup{colorlinks=true,linkcolor=blue,citecolor=cyan,urlcolor=gray}
\usepackage[nameinlink,noabbrev,capitalise]{cleveref}
\def\R{{\mathbb R}}

\def\p{\partial}
\def\O{\Omega}

\def\a{\alpha}
\def\b{\beta}

\def\d{\delta}

\def\e{\varepsilon}
\def\t{\theta}
\def\L{\Lambda}
\def\k{\kappa}

\def\SS{\mathcal{S}}

\def\l{\langle}
\def\r{\rangle}
\def\2{\mathrm{I\!I}}
\def\1{\mathrm{I}}
\def\co{\operatorname{co}}

\def\XXint#1#2#3{{\setbox0=\hbox{$#1{#2#3}{\int}$}
 \vcenter{\hbox{$#2#3$}}\kern-.5\wd0}}

\DeclareMathOperator{\supp}{supp}

\makeatletter
\newsavebox{\@brx}
\newcommand{\llangle}[1][]{\savebox{\@brx}{\(\m@th{#1\langle}\)}%
  \mathopen{\copy\@brx\kern-0.5\wd\@brx\usebox{\@brx}}}
\newcommand{\rrangle}[1][]{\savebox{\@brx}{\(\m@th{#1\rangle}\)}%
  \mathclose{\copy\@brx\kern-0.5\wd\@brx\usebox{\@brx}}}
\makeatother

\theoremstyle{plain}
\newtheorem{theorem}{Theorem}
\numberwithin{theorem}{section}
\newcommand{\newsharedstatement}[3]{%
  \newaliascnt{#1}{theorem}%
  \newtheorem{#1}[#1]{#2}%
  \aliascntresetthe{#1}%
  \crefname{#1}{#2}{#3}%
  \Crefname{#1}{#2}{#3}%
}
\newsharedstatement{lemma}{Lemma}{Lemmas}
\newsharedstatement{claim}{Claim}{Claims}
\newsharedstatement{proposition}{Proposition}{Propositions}
\newsharedstatement{cor}{Corollary}{Corollaries}
\newsharedstatement{corollary}{Corollary}{Corollaries}

\theoremstyle{definition}

\newtheorem{example}{Example}

\newsharedstatement{definition}{Definition}{Definitions}
\theoremstyle{remark}
\newtheorem{remark}{Remark}

\numberwithin{equation}{section}
\allowdisplaybreaks

\newcounter{stepctr}
{\end{list}}

\author{Artemis A. Vogiatzi}
\address{Copenhagen Centre for Geometry and Topology\\
    University of Copenhagen\\
    Universitetsparken 5, 2100 København, Danmark}
\email{aav@math.ku.dk}

\author{Qiyu Zhou}
\address{Mathematical Sciences Institute\\
Australian National University\\ Canberra ACT
2601, Australia}
\email{qiyu.zhou@anu.edu.au}

\begin{document}
\title[]{Spacelike Curve Shortening Flow of Entire Graphs in One Spacelike Codimension}
\begin{abstract}
We study curve shortening flow of entire spacelike graphs in the pseudo-Euclidean space $\R^{2,k}$, $k\geq1$, where the graph directions have mixed spacelike and timelike signature. Consequently, the region of spacelike slopes is nonconvex, making the preservation of spacelikeness more delicate. We identify a natural geometric condition that overcomes this difficulty: the closed convex hull of the initial slopes is compactly contained in the spacelike region. We prove that the convex hull is preserved along the flow, yielding uniform spacelikeness and uniform acausality. For initially $C^{1,\alpha}$ graphs, we establish global existence and uniqueness, without assuming convexity or initial curvature bounds. If the initial perturbation vanishes at spatial infinity, we further prove smooth and uniform convergence to the prescribed affine spacelike line, without rescaling. Finally, we show that every ancient entire graphical flow with globally bounded slopes that remains uniformly spacelike throughout its evolution is a stationary affine line.
\end{abstract}
\maketitle
\setcounter{tocdepth}{1}
\tableofcontents

\section{Introduction}\label{sec:introduction}

For an entire curve, global existence and convergence to a straight line are natural counterparts of the round-point behaviour of closed curves. In a pseudo-Euclidean ambient space, these questions are coupled to a further geometric problem: the evolving tangent lines must remain
spacelike. This paper gives a condition on the initial slopes that ensures global existence and uniform straightening for entire curves with one spacelike and arbitrarily many timelike normal directions. The condition allows large perturbations of an affine line and requires
no curvature convexity.

Curve shortening flow (CSF) evolves a curve $F(\cdot,t)$ according to
\begin{equation}\label{eq:intro-CSF}
    (\p_t F)^\perp = H,
\end{equation}
where $H$ is its mean curvature vector and normal projection allows for tangential reparametrisations. For closed planar curves, Gage--Hamilton \cite{GageHamilton} proved that every smooth, convex, embedded, closed plane curve shrinks to a round point under CSF. Grayson \cite{Grayson} removed the initial convexity hypothesis and showed that every smooth, embedded, closed plane curve becomes strictly convex before extinction. Andrews and Bryan~\cite{AndrewsBryan2011} later obtained 
a direct proof using a chord-distance comparison estimate. For periodic ramps in $\R^3$,
Altschuler and Grayson \cite{AltschulerGrayson1992}
established global existence and convergence to
straight lines under curve shortening flow. For free boundary curve shortening flow in high codimension, Nguyen and the first author \cite{NguyenVogiatziFBCSF} showed that if the reflected entropy is below $2$, then the flow either converges subsequentially to an orthogonal chord or develops a finite-time boundary singularity modelled on a shrinking semicircle.

For entire Euclidean graphs, the work of Ecker and Huisken \cite{EckerHuisken} established smooth, global existence for Euclidean graphical hypersurfaces with a global Lipschitz bound. The noncompact setting also requires control of the behaviour at spatial infinity to identify a limit. Miura and Rupp \cite{MiuraRupp} developed an energy method for complete Euclidean curves based on finite direction energy. Their results yield global straightening for graphical initial curves, where the positional convergence is locally smooth after translation and reparametrisation, while the derivative convergence is global in space. They also explained why uniform positional convergence need not hold in their class, see \cite[Theorem~2.2, Remark~2.3 and Corollary~2.5]{MiuraRupp}. These results distinguish the alignment of tangent directions from convergence of the whole curve to a fixed line.

Spacelike mean curvature flow has a strong global theory when all normal directions are timelike. Li and Salavessa~\cite[Theorem~1.1]{LiSalavessa} proved
long-time existence and convergence results for spacelike graphs over closed bases in pseudo-Riemannian products under some curvature assumptions. Lambert and Lotay~\cite[Theorem~1.3]{LambertLotay} established global existence for entire $n$-dimensional spacelike graphs in $\R^{n,m}$ for arbitrary timelike codimension $m$. Their theorem even allows initial graphs that approach null directions at spatial infinity and does not require
an initial mean curvature bound. 

A spacelike curve in $\R^{2,k}$ instead has normal signature $(1,k)$. The extra spacelike normal direction changes the geometry of the flow: the spacelike slope region is nonconvex, and the evolution of the induced metric contains a term of indefinite sign. Thus the entire-graph theory in $\R^{n,m}$ described above does not directly address these curves. In this mixed normal signature, Andrews and Zhou \cite{andrews2026highcodimensionmeancurvature} proved finite-time round-point convergence for compact connected spacelike-convex submanifolds. Hu and Ma \cite{HuMa} studied closed spacelike curves with a one-to-one convex projection onto a spacelike plane. They obtained a dichotomy between circular collapse and subsequential approach of tangent directions to the null cone. Our result supplies a noncompact straightening theorem in the same normal signature, under a condition on tangent slopes that requires neither spacelike convexity nor a convex projection of the curve.  

Unlike the usual spacelike graphical setting, our graphs take values in a space with both spacelike and timelike directions. We consider entire curves in $\R^{2,k}$, $k\geq 1$, written as graphs over one spacelike coordinate. The remaining graph directions have signature $(1,k)$. We fix coordinates $(x,z,y)\in \R \times \R \times \R^k$ with
\begin{equation}\label{eq:intro-metric}
    \l(x,z,y),(\widetilde{x},\widetilde{z},\widetilde{y})\r
    =x\widetilde{x} + z\widetilde{z} - y\cdot\widetilde{y}.
\end{equation}
The same fixed splitting determines the Wick-rotated norm
$\|(x,z,y)\|^2=x^2+z^2+\sum_{j=1}^k(y^j)^2$, and all norms below use this choice, which is independent of space and time. For $v,c\in\R^{1+k}$, write
\begin{equation}\label{eq:intro-graph}
    F(x,t)=(x,vx+c+\phi(x,t)),\qquad x\in\R.
\end{equation}
On the space of graph values, let
$\l p,q\r_L=p^+ q^+ -\sum_{j=1}^kp^jq^j$. The induced metric of the graph is $g(v+\phi_x)\,dx^2$, where
\begin{equation}\label{eq:intro-slope-region}
    g(p)=1+\l p,p\r_L, \qquad \SS=\{p\in\R^{1+k}:g(p)>0\}.
\end{equation}
Thus the graph is spacelike precisely when $v+\phi_x\in\SS$. In the fixed graphical parametrisation, its evolution becomes
\begin{equation}\label{eq:intro-PDE}
    \phi_t=\frac{\phi_{xx}}{g(v+\phi_x)}.
\end{equation}
The definitions and the derivation of this equation are given in the next section.

The set $\SS$ is not convex. For example, when $k=1$, the slopes $(2,2)$ and $(-2,2)$ have $g=1$, whereas their midpoint $(0,2)$ has $g=-3$. For a one-dimensional graph with only timelike graph directions, the spacelike slope region is a unit ball, whereas here the additional positive
term $(p^+)^2$ destroys this convexity. Pointwise spacelikeness therefore does not control convex combinations of tangent slopes. Since the diffusion coefficient in \cref{eq:intro-PDE} becomes singular as $g$
approaches zero, we seek a preserved region of slopes uniformly separated from the null cone. The identity in \cref{eq:g-evolution} explains the obstruction to a direct minimum-principle estimate for $g$.

The central observation is that the flow preserves the closed convex hull of the initial slopes. We write
\begin{equation}\label{eq:intro-hull}
    K_0=\overline{\co}\{v+\phi_0'(x):x\in\R\}, \qquad K_t=\overline{\co}\{v+\phi_x(x,t):x\in\R\}.
\end{equation}
Our hypothesis is $K_0\Subset\SS$. Once convex hull preservation is established, this gives a uniform positive lower bound for $g$. It also controls chords since every chord slope is an average of tangent slopes and hence belongs to $K_0$. The convex-hull condition concerns the tangent slopes and does not impose convexity on the curve itself. Consequently, with
$\k=\min_{q\in K_0}g(q) > 0$, the evolving graphs satisfy
\begin{equation}\label{eq:intro-acausal}
   \langle F(x, t)- F(y, t), F(x, t)- F (y, t)\rangle \geq \kappa(x-y)^2.
\end{equation}
This is the uniform acausality condition used below, where the same first-order hypothesis therefore controls both the local metric and the separation of distinct points on the curve. Preserving the entire convex hull retains all its supporting half-space inequalities, and provides time-independent bounds for the diffusion coefficient. This geometric estimate is the basis for the global theory.

Our main theorem is the following. The notation $C_b^{1,\a}$ denotes maps with bounded $C^{1,\a}$ norm, and $C_0$ denotes continuous maps tending to zero as $|x|\to\infty$.
\begin{theorem}\label{thm: main theorem}
Let $0 < \a <1$ and consider the initial graph $F_0(x) = (x, vx + c + \phi_0(x))$, where $v,c\in \R^{1 + k}$ and $\phi_0\in C^{1, \a}_b(\R; \R^{1 + k})$. Assume initially that $K_0 \Subset \SS$.  Then there exists a global graphical spacelike curve shortening flow  
    \begin{align*}
F(x,t)=(x,vx+c+\phi(x,t)),
    \end{align*}
where $\phi \in C([0, \infty); C^1_b(\R; \R^{1 + k})) \cap C^\infty(\R \times (0, \infty))$, $\phi(\cdot, 0) = \phi_0$. This flow is unique among classical graphical solutions with the same initial data whose perturbations are bounded and which remain uniformly spacelike on every finite time interval. Moreover, 
\begin{align*}
K_t\subseteq K_s\Subset\SS,\qquad\text{for all } 0\leq s\leq t.
\end{align*}
Assume in addition that $\phi_0\in C_0(\R; \R^{1 + k})$. Then, the flow converges smoothly and uniformly, without rescaling, to the affine spacelike line 
    \begin{equation*}
        \ell_{v,c} = \{(x, vx + c) : x\in \R\}.
    \end{equation*}
\end{theorem}
\begin{remark}
In the uniqueness assertion of Theorem~\ref{thm: main theorem}, the graphical solutions are required to have perturbations $\psi$ satisfying, on every finite time interval $[0, T]$, that
    \begin{align*}
\psi\in C(\R\times[0,T];\R^{1+k}) \cap C^{2,1}(\R\times(0,T];\R^{1+k}),
    \end{align*}
and
\begin{align*}
\sup_{\R\times[0,T]}\|\psi\|<\infty,\qquad\inf_{\R\times[0,T]}g(v+\psi_x)>0. 
\end{align*}
\end{remark}

The theorem gives a global straightening criterion for a class of spacelike curves without assuming spacelike convexity initially. Its hypothesis depends on the first derivative and permits arbitrarily large perturbations with non-proportional graph components as shown in \cref{ex:large-graphs}. The $C^{1,\a}$ formulation also permits initial curves for which curvature is not defined. 

The convergence is uniform on the whole spatial domain, including the height itself and every spatial derivative, and the affine line is fixed by the initial asymptotics, so no time-dependent translation or rescaling is needed. The assumption $\phi_0\in C_0$ specifies this positional
behaviour at the ends without imposing a spatial decay rate or an integral bound on the perturbation or its derivative. This gives a complementary convergence statement to the finite-direction-energy theory of \cite{MiuraRupp}: we prescribe a common affine asymptote at both ends and obtain uniform convergence to that line in the mixed-signature setting. 

A rigidity result holds for ancient entire graphical
flows. Unlike Theorem~\ref{thm: main theorem}, it requires
no prescribed asymptotic behaviour at spatial infinity, and assumes uniform control of the slopes and spacelikeness throughout the ancient flow.

\begin{theorem}\label{thm:ancient-rigidity}
Let $F(x,t)=(x,u(x,t))$, where $(x,t)\in\R\times(-\infty,T)$, be a smooth ancient spacelike curve shortening flow in $\R^{2,k}$, $k\geq1$. Suppose that
\begin{align}\label{eq:ancient-bounds}
\sup_{\R\times(-\infty,T)}\|u_x\|<\infty,\qquad\inf_{\R\times(-\infty,T)}g(u_x)>0.
\end{align}
Then there exist constant vectors $p_*,c_*\in\R^{1+k}$, with $g(p_*)>0$, such that
\begin{align*}
F(x,t)=(x,p_*x+c_*),
\end{align*}
for all $x\in\R$ and $t<T$.
\end{theorem}

The paper is organised as follows. Section 2 introduces the geometric setting and derives the graphical curve shortening flow equation. Section 3 establishes convex-hull preservation and quantitative acausality. In Section 4, we prove global existence and uniqueness, while Section 5 establishes smooth convergence to an affine spacelike line. Finally, Section 6 proves the ancient graphical rigidity theorem.

\section*{Declarations} 
\subsection*{Acknowledgments.} The first author is financially supported by the Danish National Research Foundation through the Copenhagen Centre for Geometry and Topology (DNRF151). 

\subsection*{AI disclosure.} The authors used LLMs to help with checking mathematical arguments, identifying relevant literature and revising the exposition. The authors take full responsibility for the mathematical content, references and final manuscript.

\section{Preliminaries}

\subsection{Pseudo-Euclidean space and spacelike curves}
For integers $n\geq 1$ and $k\geq 0$, the pseudo-Euclidean space $\R^{n,k}$ is $\R^{n+k}$ equipped with the nondegenerate bilinear form
    \begin{equation*}
        \l X,Y\r =\sum_{i=1}^{n}X^iY^i -\sum_{j=1}^{k}X^{n+j}Y^{n+j}.
    \end{equation*}
Its Levi-Civita connection is the flat connection in these coordinates. A nonzero vector $X$ is called spacelike, timelike or null if  $\l X,X\r$ is positive, negative or zero respectively. An immersion is spacelike if the restriction of this bilinear form to each tangent space is positive definite.

Throughout this paper, the ambient space is $\R^{2,k}\cong \R\times\R^{1,k}$, with $k\geq 1$. We write a point as $(x,p)$ where $p=(p^+,p^1,\ldots,p^k)$, and use
    \begin{equation*}
        \l p,q\r_L=p^+q^+-\sum_{j=1}^k p^jq^j, \qquad \l(x,p),(y,q)\r = xy+\l p,q\r_L.
    \end{equation*}
Thus $\R^{1+k}$ denotes the space of graph values in coordinates, and $\l\cdot,\cdot\r_L$ specifies its metric of signature $(1,k)$. The fixed coordinate splitting determines the Wick-rotated positive definite inner product
    \begin{equation*}
        \llangle(x,p),(y,q)\rrangle =xy+p^+q^++\sum_{j=1}^k p^jq^j,
    \end{equation*}
and its norm is denoted by $\|\cdot\|$. On the space of graph values,
    \begin{equation*}
        \|p\|^2 = (p^+)^2 + \sum_{j=1}^k(p^j)^2.
    \end{equation*}
We keep this inner product fixed throughout the flow. All norms of graph-valued maps below, including supremum and H\"older norms, are with respect to this fixed Wick-rotated norm. For linear functionals and linear maps, $\|\cdot\|$ denotes the corresponding dual and operator norms.

Let $I\subset\R$ be an interval and let $F:I \to\R^{2,k}$ be a smooth parametrised curve. It is spacelike if $\l F_x,F_x\r > 0 $ at every point. Its arclength parameter $s$ is determined by 
    \begin{equation*}
        ds= \sqrt{\l F_x,F_x\r}\,dx, \qquad \p_s=\frac{1}{\sqrt{\l F_x,F_x\r}} \p_x, 
    \end{equation*}
where the arclength is computed using the induced pseudo-Euclidean metric. We define the mean curvature vector by
    \begin{equation*}
        H:= \p_s^2F =\frac{1}{\sqrt{\l F_x,F_x\r}}\p_x \left(\frac{F_x}{\sqrt{\l F_x,F_x\r}}\right), 
    \end{equation*}
which is normal to the curve since
\begin{align*}
    \l H,\p_sF\r=\frac{1}{2}\p_s\l\p_sF,\p_sF\r=0,
\end{align*}
where $\l\p_sF,\p_sF\r=1$.

The normal space has signature $(1,k)$, and the orthogonal projection of an ambient vector $V$ onto it is
    \begin{equation*}
        V^\perp =V-\frac{\l V,F_x\r}{\l F_x,F_x\r}F_x.
    \end{equation*}
For a smooth family $F(x,t)$, these definitions apply at each fixed time, and the family moves by CSF if
    \begin{equation*}
        (\p_tF)^\perp = H.
    \end{equation*}

\subsection{One-dimensional graphs and graphical CSF}
Let $u=(u^+,u^1,\ldots,u^k)$ be a graph function which is smooth for positive times, and write $F(x,t)=(x,u(x,t))$ for $x\in \R$. The coordinate derivatives are
    \begin{equation*}
        F_x=(1,u_x),\qquad F_{xx}=(0,u_{xx}),\qquad F_t=(0,u_t).
    \end{equation*}
The induced metric is
    \begin{equation}
        F(\cdot,t)^*\l\cdot,\cdot\r
        =g(u_x)\,dx^2,\qquad g(p)=1+\l p,p\r_L =1+(p^+)^2-\sum_{j=1}^k(p^j)^2, 
    \end{equation}
where the graph is spacelike when $g(u_x)>0$. The arclength satisfies
    \begin{equation*}
        ds= \sqrt{g(u_x)}\,dx,\qquad \p_s=g(u_x)^{-1/2}\p_x.
    \end{equation*}
    We use
    \begin{equation*}
        \SS=\{p\in\R^{1+k}:g(p)>0\},
    \end{equation*}
for the open set of spacelike slopes. A graphical family is \emph{uniformly} spacelike on a time interval $I$ if $\inf_{\R\times I}g(u_x) > 0$. We call a graph \emph{acausal} if every ambient chord joining distinct points is strictly spacelike, and \emph{uniformly} acausal if there exists some $\k > 0$ such that 
\begin{equation*}
    \l F(x, t)-F(y, t), F(x, t)-F(y, t)\r \geq \kappa(x-y)^2,
\end{equation*}
for all $x,y\in \R$, $t\in I$. 

In the following calculation we abbreviate $g=g(u_x)$. In this notation, $\l F_x,F_x\r=g$. Using $g_x=2\l u_x,u_{xx}\r_L$, the mean curvature of the graph is 
    \begin{equation}\label{eq:prelim-mean-curvature}
        \begin{aligned}
            H &=\frac{F_{xx}}{g}-\frac{g_x}{2g^2}F_x=\left(0,\frac{u_{xx}}{g}\right) -\frac{\l u_x,u_{xx}\r_L}{g^2}(1,u_x)=\left(0, \frac{u_{xx}}{g} \right)^\perp.
        \end{aligned}
    \end{equation}
The CSF equation is then
    \begin{equation*}
        \left(0,u_t-\frac{u_{xx}}{g}\right)^\perp = 0.
    \end{equation*}
Hence $(0,u_t-u_{xx}/g)$ is tangent to $F(\R,t)$ and equals $\lambda(1,u_x)$ for some scalar $\lambda \in \R$. Comparing the first coordinates gives $\lambda=0$, and then the remaining coordinates give
    \begin{equation}\label{graphical MCF equation}
        u_t=\frac{u_{xx}}{g(u_x)}=\frac{u_{xx}}{1+(u_x^+)^2- \sum_{j=1}^k(u_x^j)^2}.
    \end{equation}
Its velocity in the fixed graphical parametrisation is
    \begin{equation*}
        F_t= H + \frac{\l u_x,u_{xx}\r_L}{g^2}F_x.
    \end{equation*}
The last term is tangential and accounts for the choice of the fixed coordinate $x$. Each of the $1+k$ components of $u$ satisfies the same scalar diffusion equation with coefficient $g(u_x)^{-1}$.

For the corresponding graphical formulation in higher codimension with timelike graph directions, see \cite[Appendix~A]{LambertLotay}.

\subsection{Affine reference lines and function spaces}
For fixed $v,c\in\R^{1+k}$, let
    \begin{equation*}
        \ell_{v,c}=\{(x,vx+c):x\in\R\}.
    \end{equation*}
Its tangent vector is $(1,v)$, so this line is spacelike if and only if $g(v) > 0$. Writing $u(x,t)= vx+c+\phi(x,t)$, where $\phi\colon \mathbb R\times[0,\infty)\to\mathbb R^{1+k}$, we have 
    \begin{equation*}
        u_t=\phi_t,\qquad u_x=v+\phi_x,\qquad u_{xx}= \phi_{xx}.
    \end{equation*}
Substituting this into \cref{graphical MCF equation} gives
    \begin{equation}\label{eq: prelim-perturbation}
        \phi_t =\frac{\phi_{xx}}
        {1+(v^+ + \phi_x^+)^2 - \sum_{j=1}^k(v^j+\phi_x^j)^2}, \qquad \phi(\cdot,0) = \phi_0,
    \end{equation}
where
\begin{align*}
    \l v+\p_x\phi,v+\p_x\phi\r_L=\left(v^{+}+\p_x\phi^{+}\right)^2-\sum_{j=1}^k \left(v^j+\p_x\phi^j\right)^2. 
\end{align*}
We write $p= v + \phi_x$ for the slope of the full graph and $a(p)=1/g(p)$ on $\SS$. Thus \cref{eq: prelim-perturbation} is $\phi_t=a(v+\phi_x)\phi_{xx}$.

For $m=1,2$ and $0<\alpha<1$, the notation
$C_b^{m+\alpha}(\R;\R^{1+k})$ denotes the space of $C^m$ maps whose derivatives through order $m$ are bounded and whose derivatives of order $m$ have finite global $\alpha$-H\"older seminorm. We write
    \begin{equation*}
        [h]_\a =\sup_{x\neq y}\frac{\|h(x)-h(y)\|}{|x-y|^\a}, \qquad \|h\|_\infty=\sup_{x\in\R}\|h(x)\|.
    \end{equation*}
Here $\|\cdot\|$ is the fixed Wick-rotated norm defined before. The space $C_0(\R;\R^{1+k})$ consists of
continuous maps tending to zero as $|x|\to\infty$. The notation $C^{2,1}$ means that the spatial derivatives through order two and the first time derivative are continuous on the indicated domain.

The closed convex hull of the slopes at time $t$ is
    \begin{equation*}
        K_t=\overline{\co}\{v+\phi_x(x,t):x\in\R\}, \qquad
        K_0=\overline{\co}\{v+\phi_0'(x):x\in\R\}.
    \end{equation*}
The convex hull assumption does not require the initial curve to be convex or close to an affine line. Rather, it imposes a condition on the slopes of the graph, ensuring that every convex combination of these slopes remains uniformly spacelike. In particular, it allows initial curves that are arbitrarily large perturbations of a straight line, involving both spacelike and timelike graph directions, and whose curvature varies along the curve. To illustrate this, we give the following example.

\begin{example}\label{ex:large-graphs}
Consider in $\R^{2,1}$ the initial graph
    \begin{align*}
F_0(x)=\bigl(x,2Ae^{-x^2},Ae^{-x^2}+\e e^{-4x^2}\bigr),
\qquad A>0,\quad 0<\e<\sqrt{\frac{3e}{32}}.
    \end{align*}
The parameter $A$ can be arbitrarily large, so the perturbation is not required to be small. Moreover, the Gaussian profile changes concavity, and the additional term $\e e^{-4x^2}$ makes the two graph components nonproportional. Thus, this example is not restricted to curves with a fixed curvature sign or to the special case in which the contributions to the induced metric cancel.

Nevertheless, the convex-hull condition is satisfied. Indeed, let
\begin{align*}
s(x)=(e^{-x^2})',\qquad r(x)=(e^{-4x^2})'.
\end{align*}
Every slope in $K_0$ has the form $(2As,As+\e r)$, where
\begin{align*}
(s,r)\in\overline{\operatorname{co}}\left\{(s(x),r(x)):x\in\R\right\}.
\end{align*}
Since $|r(x)|\leq\sqrt{8/e}$ for every $x$, the same bound holds for every such $r$. Consequently,
\begin{align*}
g(2As,As+\e r)
&=1+4A^2s^2-(As+\e r)^2\\
&=1+3\left(As-\frac{\e r}{3}\right)^2
-\frac{4}{3}\e^2r^2\\
&\geq 1-\frac{32\e^2}{3e}>0.
\end{align*}
Importantly, the lower bound is independent of $A$, and therefore remains valid even for arbitrarily large perturbations. This shows that the convex-hull condition controls the spacelike geometry of the curve without imposing smallness or convexity of the initial graph.
\end{example}

\begin{example} 
In $\R^{2,1}$ take the initial curve $F_0(x)=(x,e^{-x^2},e^{-x^2})$. Its solution to CSF is
    \begin{equation*}
        F(x,t)=(x,h(x,t),h(x,t)),\qquad h(x,t)=\frac{1}{\sqrt{1+4t}} \exp\left(-\frac{x^2}{1+4t} \right).
    \end{equation*}
The induced metric is $g=1+h_x^2-h_x^2 = 1$, so the graphical flow equation reduces to $h_t=h_{xx}$.
    
Every vertical slope is $(h_x,h_x)$. At $t=0$, $h_x= -2x e^{-x^2}$ has range $[-\sqrt{2/e},\sqrt{2/e}]$, so the initial convex hull of slopes is the line segment
    \begin{equation*}
        K_0=\bigl\{(s,s): |s|\leq\sqrt{2/e}\bigr\}\Subset \SS.
    \end{equation*}
At time $t$, the same calculation gives
    \begin{equation*}
        K_t=\left\{(s,s): |s|\leq\frac{\sqrt{2/e}}{1+4t}\right\}.
    \end{equation*}
Thus $K_t$ shrinks to $\{(0,0)\}$, and since $\|h(\cdot,t)\|_\infty=(1+4t)^{-1/2}$, the curve converges smoothly and uniformly to the spacelike axis $F_\infty(x)=(x,0,0)$.
\end{example}

\section{Graphical flow and convex hull preservation}
Our first goal is to show that the convex hull of the initial slopes is preserved by the flow. This is important because the spacelike slope region is not convex, and we need to prevent the evolving slopes from approaching the null cone. 

The evolution of the induced metric explains why a direct
minimum-principle argument does not suffice. For a smooth solution, write $p = v+\phi_x$, $ a = a(p) =1/g(p)$ and $g=g(p)$. Differentiating the graphical equation gives $p_t = \p_x(a p_x)$, while
\begin{equation*}
    g_x=2\l p,p_x\r_L,\qquad a_x=-a^2 g_x.
\end{equation*}
Consequently,
\begin{equation}\label{eq:g-evolution}
    \begin{aligned}
    (\p_t-a\p_{xx})g
      &=2a_x\l p,p_x\r_L-2a\l p_x,p_x\r_L\\
      &=-a^2(g_x)^2-2a\l p_x,p_x\r_L.
    \end{aligned}
\end{equation}
The last inner product has no fixed sign. Thus, rather than estimating the induced metric directly, we will use the scalar nature of the graphical equation to control the slopes. The following lemma is a maximum principle on unbounded domains, which will be used throughout the paper.
\begin{lemma}\label{lem : maximum principle for slope preservation}
Let $\O \subset \R^d$ be open and let $f\in C(\overline{\O} \times[0,T])\cap C^{2,1}(\O \times(0,T])$
be bounded above. Suppose
    \begin{equation*}
        f_t-b^{ij}(z,t)f_{z_i z_j}\leq Bf, \qquad 0\le (b^{ij})\leq \L I,
    \end{equation*}
where the matrix is symmetric and $B\geq 0$ is constant. If $f\leq 0$ initially and on the lateral boundary, then $f\leq 0$ throughout the cylinder.
\end{lemma}
\begin{proof}
Replacing $f$ by $e^{-Bt}f$ reduces to $B=0$. Let
    \begin{equation*}
        \Psi(z,t)=e^{(2d\L +1)t}(1+|z|^2). 
    \end{equation*}
We first note that the spatial Hessian of $\Psi$ is $2e^{(2d\L+1)t}I$. Since $(b^{ij})\leq \L I$, each diagonal entry satisfies $b^{ii} \leq \L$, and hence $\sum_{i=1}^d b^{ii}\le d\L$. Differentiating in time and using this trace bound gives
    \begin{align*}
        e^{-(2d\L+1)t}\left(\p_t-b^{ij}\p_{z_i z_j}\right) \Psi & = (2d\L+1)(1+|z|^2)-2\sum_{i=1}^d b^{ii}\\
        &\ge 1+(2d\L+1)|z|^2\geq 1.
    \end{align*}
For any $\e > 0$ let $w = f-\e \Psi$. The initial
and lateral boundary values of $w$ are strictly negative. Since $f$ is bounded above and $\Psi \geq 1+|z|^2$, there is some $R <\infty$ such that $w<0$ whenever $|z|\geq R$, uniformly over $[0,T]$. If $w$ were positive somewhere, it would attain a positive maximum on the compact set $(\overline{\O} \cap\overline{B_R(0)}) \times [0,T]$. At such points we have $(\p_t-b^{ij}\p_{z_i z_j}) w \geq 0$, which is a contradiction because at every interior point we have 
    \begin{equation*}
        (\p_t-b^{ij}\p_{z_i z_j}) w = (\p_t-b^{ij}\p_{z_i z_j})f - \e (\p_t-b^{ij}\p_{z_i z_j}) \Psi \leq -\e e^{(2d\L + 1)t} < 0. 
    \end{equation*}
Thus $f\leq \e \Psi$, letting $\e \to 0$ at each point gives the result.
\end{proof}

We now apply the maximum principle to the graphical flow. The key observation is that every component of the graph satisfies the same scalar parabolic equation. By considering differences between two points on the graph, we can control their secant slopes and consequently the tangent slopes, rather than attempting to apply a maximum principle to the indefinite quantity $g(p)$. This gives the following preservation result.

\begin{proposition}[Convex hull preservation]\label{prop: hull preservation}
Suppose $u=vx+c+\phi$, where $\phi\in C(\R\times[0,T])\cap C^{2,1}(\R\times(0,T])$ is bounded, $u(\cdot,0)$ is continuously differentiable with bounded derivative, and $u_t=b(x,t)u_{xx}$ for $0 < t\leq T$. Assume $0 < b\leq \Lambda_T$. If $K_0$ is the closed convex hull of the initial slopes, then
    \begin{equation*}
        \frac{u(x,t)-u(y,t)}{x-y}\in K_0\quad( x \neq y), \qquad u_x(x,t)\in K_0\quad(0\leq t\leq T).
    \end{equation*}
\end{proposition}
\begin{proof}
First we show $v\in K_0$. Taking an average we have
    \begin{equation*}
        \frac{1}{2R} \int_{-R}^R u_x(x,0)\,dx =v + \frac{\phi_0(R)-\phi_0(-R)}{2R} \longrightarrow v, \qquad \text{as } R \longrightarrow \infty, 
    \end{equation*}
where the average over $[-R,R]$ is the limit of Riemann sums with nonnegative weights, each sum belongs to $K_0$, hence $v\in K_0$ follows from the closedness of $K_0$ and that $\phi_0$ is bounded. 
    
Let $L : \R^{1+k}\to\R$ be a linear functional and
$h=\max_{q\in K_0}L(q)$. On the open half-plane $x > y$, define
    \begin{equation*}
        Z(x,y,t)= L(u(x,t) - u(y,t)) - h(x-y).
    \end{equation*}
We note that $Z=0$ on the diagonal. On the initial slice the fundamental theorem of calculus gives
    \begin{equation*}
        Z(x,y,0)= \int_y^x [L(u_x(z,0))-h] \,dz \leq 0.
    \end{equation*}
For positive time, we have
    \begin{equation*}
        \begin{aligned}
            Z(x,y,t) &=L\bigl(u(x,t)-u(y,t)\bigr)-h(x-y)\\
            &=(x-y)L(v)+L\bigl(\phi(x,t)-\phi(y,t)\bigr)-h(x-y)\\
            &=(L(v)-h)(x-y)+L\bigl(\phi(x,t)-\phi(y,t)\bigr).
        \end{aligned}
    \end{equation*}
Here $L(v)\leq h$, $x - y > 0$, and the second term is bounded on the slab, so $Z$ is bounded above. Moreover,
\begin{align*}
Z_t&=L\bigl(b(x,t)u_{xx}(x,t)-b(y,t)u_{xx}(y,t)\bigr),\\
Z_{xx}&=L(u_{xx}(x,t)),\\
Z_{yy}&=-L(u_{xx}(y,t)).
\end{align*}
Consequently,
\begin{align*}
Z_t-b(x,t)Z_{xx}-b(y,t)Z_{yy}=0.
\end{align*}
The coefficient matrix is $\operatorname{diag}(b(x,t),b(y,t))$, which satisfies the hypotheses of \cref{lem : maximum principle for slope preservation}. Hence, \cref{lem : maximum principle for slope preservation} gives $Z \leq 0$.
    
Dividing by $x - y>0$ shows that the secant slope belongs to the half-space $\{ q : L(q) \leq h\}$, and the result holds for every $L$. By the finite-dimensional separating hyperplane theorem\footnote{which says if $p$ lies outside a nonempty closed convex set in a finite dimensional vector space, then there exists a linear functional such that $L(p) > \sup_{q\in K_0} L(q)$.} \cite[§2.5.1, Example 2.20, p. 49]{BoydVandenberghe}, a point outside $K_0$ violates at least one of these inequalities, so the secant slope lies in $K_0$. Reversing the two endpoints gives the result for $x < y$. Taking $x\to y$ and using closedness proves the claim for the slopes. 
\end{proof}

Applying the proposition with time $s$ as the initial time gives
\begin{align*}
K_t\subseteq K_s,\qquad 0\leq s\leq t.
\end{align*}
With $K_0\Subset\SS$, continuity and positivity of $g$ give $0 < \k= \min_{K_0}g \leq \max_{K_0} g = G < \infty$. In particular,
\begin{equation}\label{rangeofa}
    \lambda :=G^{-1}\leq a(p(x,t))\leq \k^{-1}=: \L.
\end{equation}
All these constants depend only on the initial convex hull. Since every component of $\phi$ satisfies
$\phi_t^i=a(p)\phi_{xx}^i$, we have
\begin{align*}
(\p_t-a(p)\p_{xx})\|\phi\|^2&=2\sum_i\phi^i\bigl(\phi_t^i-a(p)\phi_{xx}^i\bigr)-2a(p)\sum_i|\phi_x^i|^2\\
&=-2a(p)\|\phi_x\|^2\leq0.
\end{align*}
Thus the function
\begin{align*}
f(x,t)=\|\phi(x,t)\|^2-\|\phi_0\|_\infty^2
\end{align*}
satisfies
\begin{align*}
(\p_t-a(p)\p_{xx})f\leq0,\qquad f(\cdot,0)\leq0.
\end{align*}
Applying \cref{lem : maximum principle for slope preservation} on each finite time interval gives
\begin{align*}
\|\phi(\cdot,t)\|_\infty\leq\|\phi_0\|_\infty.
\end{align*}
Convex hull preservation also gives a geometric consequence. Since every chord slope is an average of tangent slopes, it remains in $K_0$. Therefore, not only are the tangent vectors spacelike, but any two distinct points on the evolving curve are spacelike separated. The following proposition makes this quantitative.

\begin{proposition}\label{prop: acausality}
Let $F(x,t)=(x,u(x,t))$ be a spacelike graph. If $u_x(\cdot,t)$ takes values in a closed convex set $K_0\Subset\SS$, then
    \begin{equation*}
      \langle F(y, t)-F(x, t), F(y, t)-F(x, t)\rangle \geq \kappa(y-x)^2,
    \end{equation*}
    where $\k = \min_{K_0}g$. 
\end{proposition}
\begin{proof}
For $x\neq y$, the secant slope satisfies
    \begin{equation*}
        q(x,y,t):=\frac{u(y,t)-u(x,t)}{y-x} =\int_0^1 u_x(x+\t (y-x),t) \,d\t\in K_0
    \end{equation*}
by closedness and convexity. Consequently,
    \begin{equation*}
    \langle F(y, t)-F(x, t), F(y, t)-F(x, t)\rangle=(y-x)^2 g(q(x, y, t)) \geq \kappa(y-x)^2.
    \end{equation*}
The diagonal case is immediate.
\end{proof}

\section{Existence and uniqueness}
For noninteger $s>0$, let $h_b^s$ be the closure in $C_b^s$ of the smooth maps with all derivatives bounded. For the little H\"older spaces used below, interpolation simply gives an intermediate order of spatial regularity, the precise identity needed in the proof is provided by \cite[Proposition~2.2(c)]{ShaoSimonett}.

\begin{lemma}[Local existence for $C^{1,\alpha}$ initial data]\label{lem: local existence}
Let $\phi_0\in C_b^{1,\a}(\R^n;\R^m)$, let $V\in\R^{m\times n}$ and $B(P)$ be a smooth symmetric matrix function on $\R^{m\times n}$ with bounded derivatives and $\lambda I \leq B(P)\leq \L I$ for constants $0< \lambda \leq \L<\infty$. Fix
$0<\sigma<\b<\a$ and put $\mu=(1+\b-\sigma)/2$. The system
    \begin{equation*}
     \p_t\phi=B^{ij}(V+D\phi)\p_{ij}\phi, \qquad \phi(\cdot,0) =\phi_0,
    \end{equation*}
has a local solution which is smooth for positive times and belongs to $C([0,T];h_b^{1+\b}) \cap C((0,T];h_b^{2+\sigma}) \cap C^1((0,T];h_b^\sigma)$ for some $T > 0$. Moreover, 
    \begin{equation}\label{eq: local-weight}
        \sup_{0<t\le T} t^{1-\mu} \big(\|\phi(t)\|_{C_b^{2+\sigma}}+ \|\p_t\phi(t)\|_{C_b^\sigma}\big)<\infty.
    \end{equation}
The solution continues through any finite endpoint at which it has a limit in $h_b^{1+\b}$.
\end{lemma}
\begin{proof}
Let $E_0 = h_b^\sigma$, $E_1 = h_b^{2 + \sigma}$ and $E_\mu = (E_0,E_1)^0_{\mu,\infty}$. The inequality $0<\sigma<\b<1$ gives $1/2<\mu<1$, and
$\sigma+2\mu=1+\b$. Consequently, the interpolation identity in \cite[Proposition~2.2(c)]{ShaoSimonett} gives
$E_\mu = h_b^{1+\b}$ with equivalent norms. The loss from $\a$ to $\b$ ensures that the initial value belongs to this space, since 
    \begin{equation*}
        \|\phi_{0,\e}-\phi_0\|_{C_b^{1+\b}} \leq C\e^{\a-\b}\|\phi_0\|_{C_b^{1,\a}} \longrightarrow 0, 
    \end{equation*}
where $\phi_{0,\e}$ is the standard mollification of $\phi_0$. 

For $w\in E_\mu$ and $h\in E_1$, define $A(w)h=-B^{ij}(V+D w)\p_{ij}h$. Since $Dw\in h_b^\b$ and $B$ is smooth, each coefficient $B^{ij}(V+Dw)$ belongs to $h_b^\b\subset h_b^\sigma$. The H\"older product estimate $\|fg\|_{C_b^\sigma} \leq \|f\|_{C_b^\sigma}\|g\|_{C_b^\sigma}$ therefore shows that
$A(w):E_1\to E_0$ is bounded. Differentiating the composition with $B$ and applying the same product
estimate at each order proves that $w \mapsto A(w)$ is smooth.

There is a separate linear-theory hypothesis to check: the scalar operator acting on each component is uniformly elliptic, so the continuous maximal-regularity theorem \cite[Theorem~3.6]{ShaoSimonett} applied componentwise on Euclidean space shows that $A(w)$ has continuous maximal regularity for the pair $(E_0,E_1)$ and the time weight $t^{1-\mu}$. This verifies the linear existence hypothesis of the quasilinear theorem.

We may now apply \cite[Theorem~3.1(a)]{ClementSimonett} to
    \begin{equation*}
        \p_t\phi+A(\phi)\phi=0, \qquad \phi(0)=\phi_0,
    \end{equation*}
with the open phase space equal to all of $E_\mu$. It gives a solution in $C([0,T];h_b^{1+\b}) \cap C((0,T];h_b^{2+\sigma}) \cap C^1((0,T];h_b^\sigma)$, and in fact gives the stronger weighted condition
    \begin{equation}\label{eq: local-weight-vanishing}
        \lim_{t\to 0}t^{1-\mu} \bigl(\|\phi(t)\|_{E_1}+\|\p_t\phi(t)\|_{E_0}\bigr) = 0.
    \end{equation}
Continuity away from zero and \cref{eq: local-weight-vanishing} imply \cref{eq: local-weight}. 

Smooth dependence of $A$ now allows us to apply
\cite[Theorem~6.1]{ClementSimonett}, which gives
$\phi\in C^\infty((0,T);E_1)$. On any closed time interval away from zero, $D^2\phi$ and
$B(V+D\phi)$ therefore have parabolic $\sigma$-H\"older regularity. The equation has the same scalar principal part in every component, thus standard interior regularity for uniformly parabolic quasilinear equations \cite[Ch.~III, \S12, Theorem~12.2 and Ch.~IV, \S10,
Theorem~10.1]{LSU} now implies that $\phi$ is smooth for positive times.

Maximal continuation and the stated extension criterion follow from \cite[Theorem~4.1(a),(c)]{ClementSimonett}, since the phase space is all of $E_\mu$.
\end{proof}

This gives the short time existence of the flow.  

\begin{proposition}[Short time existence]\label{prop: short time existence}
Under the assumptions of \cref{thm: main theorem}, there is a $T_0>0$ and a perturbation
    \begin{equation*}
        \phi\in C([0,T_0];C_b^1(\R;\R^{1+k})) \cap C^\infty(\R\times(0,T_0);\R^{1+k})
    \end{equation*}
with initial value $\phi_0$ such that
    \begin{equation*}
        \phi_t=\frac{\phi_{xx}}{g(v+\phi_x)} \qquad\text{on } \R\times(0,T_0).
    \end{equation*}
The graph $F(x,t)=(x,vx+c+\phi(x,t))$ is uniformly spacelike on $[0,T_0]$. 
\end{proposition}
\begin{proof}
We first extend the coefficient to all slopes so that \cref{lem: local existence} applies, then continuity at the initial time recovers the original equation. For $r>0$, write
    \begin{equation*}
        \mathcal{N}_r =\{p\in\R^{1+k}: d(p,K_0)\leq r\},
    \end{equation*}
where distance is measured using the fixed Wick-rotated norm. Since $K_0$ is compact and contained in the open set $\SS$, we can choose $\rho > 0$ such that $\mathcal{N}_{2\rho}\Subset\SS$. Choose $\chi\in C_c^\infty(\SS)$ with $0\leq \chi \leq 1$ and $\chi=1$ on $\mathcal{N}_\rho$. Extend $\chi$ and $\chi a$ by zero outside $\SS$, and set $\widetilde{a}=\chi a+(1-\chi)$ on $\R^{1+k}$. Since $0\leq\chi\leq1$ and $a$ is smooth
and strictly positive on $\supp\chi$, the extension $\widetilde a$ is smooth with bounded derivatives and satisfies
\begin{align*}
0<\widetilde\lambda\leq\widetilde a(p)\leq\widetilde\Lambda<\infty,\qquad p\in\R^{1+k}.
\end{align*}
Consider
    \begin{equation}\label{eqwithtildea}
        \phi_t=\widetilde{a} (v+\phi_x)\phi_{xx}, \qquad \phi(\cdot,0)=\phi_0.
    \end{equation}
Fix $0 < \sigma<\b <\a$ and $\mu=(1+\b-\sigma)/2$. Apply \cref{lem: local existence} with $n=1$, $m=1+k$, $V=v$ and $B^{11}=\widetilde{a}$. Denote the maximal interval of this solution by $[0,T_{\max})$, the solution is smooth for positive times and continuous into $h_b^{1+\b}$ at $t=0$. Since this space embeds continuously into $C_b^1$, we obtain
    \begin{equation*}
        \|\phi(\cdot,t)-\phi_0\|_\infty + \|\phi_x(\cdot,t)-\phi_0'\|_\infty \longrightarrow 0,
        \qquad \text{as } t \to 0.
    \end{equation*}
Choose $0<T_0<T_{\max}$ sufficiently small that
    \begin{equation*}
        \sup_{0\leq t\leq T_0} \|\phi_x(\cdot,t)-\phi_0'\|_\infty < \rho.
    \end{equation*}
For every $x\in\R$, the initial slope $v+\phi_0'(x)$ belongs to $K_0$. Hence, for $0\leq t\leq T_0$,
    \begin{equation*}
        d(v+\phi_x(x,t),K_0) \leq \|\phi_x(x,t)-\phi_0'(x)\|< \rho.
    \end{equation*}
Every evolving slope therefore belongs to $\mathcal{N}_\rho$, where $\widetilde{a} = a = 1/g$. Compactness of $\mathcal{N}_\rho$ also gives the quantitative bound
    \begin{equation*}
        g(v+\phi_x(x,t)) \geq \min_{p\in \mathcal{N}_\rho}g(p)>0 \qquad x\in\R,\ 0\leq t\leq T_0.
    \end{equation*}
This proves uniform spacelikeness, while continuity into $C_b^1$ and positive-time smoothness give the regularity.
\end{proof}

Next we prove higher derivative estimates. 

\begin{proposition}[Uniform a priori estimates]
\label{prop: global estimates}
Let $\phi$ be the maximal solution obtained from \cref{prop: short time existence}, defined on $[0,T_{\max})$. Then
\begin{align}\label{phi0infty}
\|\phi(\cdot,t)\|_\infty\leq\|\phi_0\|_\infty.
\end{align}
For every $\tau>0$ and integer $m\geq0$, with $\tau<T_{\max}$, there exists $C_{m,\tau}<\infty$ such that
\begin{align}\label{derest}
\sup_{\tau\leq t<T_{\max}}\|\p_x^m\phi(\cdot,t)\|_\infty\leq C_{m,\tau}.
\end{align}
\end{proposition}
\begin{proof}
We first establish uniform estimates on the maximal existence interval. For every $T<T_{\max}$, continuity of
$\phi$ into $h_b^{1+\beta}$ on $[0,T]$ implies that $\phi$ is bounded on $\R\times[0,T]$. We may therefore apply \cref{prop: hull preservation} to the solution of the extended \cref{eqwithtildea}, obtaining
\begin{align*}
v+\phi_x(x,t)\in K_0,\qquad x\in\R,\quad 0\leq t\leq T.
\end{align*}
Since $T<T_{\max}$ was arbitrary, this holds for every $0\leq t<T_{\max}$. Since $\widetilde a=a$ on $K_0$, the solution of the modified equation also solves
\begin{align}\label{originaleq}
\begin{cases}
\phi_t=a(v+\phi_x)\phi_{xx},\\ 
\phi(\cdot,0)=\phi_0.
\end{cases}
\end{align}
Moreover, since $K_0\Subset\mathcal S$ and $g$ is continuous and strictly positive on $\mathcal S$, there exist constants $0<c_1\leq c_2<\infty$, depending only on $K_0$, such that
\begin{align*}
c_1\leq g(v+\phi_x(x,t))\leq c_2.
\end{align*}
Consequently,
\begin{align*}
0<\lambda:=\frac1{c_2}\leq a(v+\phi_x(x,t))\leq\frac1{c_1}=:\Lambda<\infty.
\end{align*}
Furthermore, using the fixed Wick-rotated norm on $\R^{1+k}$, we compute
\begin{align*}
(\p_t-a\p_{xx})\|\phi\|^2=-2a\|\phi_x\|^2\leq0
\end{align*}
and since $0<\lambda\leq a\leq\Lambda$, the scalar parabolic maximum principle on $\R$ gives \cref{phi0infty}, which means that the solution remains uniformly bounded for all times for which it exists.

Recall that
\begin{align*}
p_t=\partial_x(a(p)p_x),\qquad \lambda\leq a(p)\leq\Lambda,
\end{align*}
and that $p$ is uniformly bounded since $p\in K_0$. Each component of $p$ therefore satisfies a uniformly parabolic scalar divergence-form equation.

Fix $0<\tau<T_{\max}$ and choose $r>0$ sufficiently small that $16r^2<\tau/2$. For every $t\in[\tau,T_{\max})$, the backward parabolic cylinder
\begin{align*}
Q_{4r}^-(x,t)=(x-4r,x+4r)\times(t-16r^2,t]
\end{align*}
lies in $\R\times(\tau/2,T_{\max})$. Applying Nash's interior Hölder estimate on these cylinders gives an exponent $\delta\in(0,1)$ and uniform local
$C^{\delta,\delta/2}$ bounds for $p$. The constants depend on $\tau$, the ellipticity bounds and the initial data, but not on the centre $(x,t)$ or on an upper time endpoint $S<T_{\max}$.

Since $a$ is smooth on a neighbourhood of $K_0$, the same Hölder bounds hold for $a(p)$. Applying interior parabolic Schauder estimates on smaller cylinders therefore yields
\begin{align}\label{c2estimate}
\|\phi\|_{C^{2+\delta,1+\delta/2}(\R\times[\tau,S])}\leq C_\tau,
\end{align}
for every $\tau<S<T_{\max}$, with $C_\tau$ independent of $S$.

From \cref{c2estimate}, we get
\begin{align}\label{boundedphixx}
\sup _{\tau \leq t<T_{\max }}\left\|\phi_{x x}(\cdot, t)\right\|_{\infty} \leq C_\tau.
\end{align}
To estimate $\phi_{xxx}$, we return to the equation for $p=v+\phi_x$. The previous estimate controls $p_x=\phi_{xx}$, allowing us to apply Schauder theory once more. We write
\begin{align*}
p_t=a(p) p_{xx}+\left(D a(p) \cdot p_x\right) p_x .
\end{align*}
We already control $p_x=\phi_{xx}$ in Hölder norm. Hence the coefficients and the last term are Hölder continuous. Consequently, applying \cref{c2estimate} with $\tau$ replaced by $\tau/2$, we control $p_x=\phi_{xx}$ in parabolic Hölder norm on $[\tau/2,S]$,
\begin{align*}
\|\left(D a(p) \cdot p_x\right) p_x\|_{C^{\d, \d / 2}(\mathbb{R} \times[\tau/2,S])} \leq C_\tau.
\end{align*}
We can therefore apply the Schauder estimate once more, now to the equation for $p$, obtaining
\begin{align*}
\|p\|_{C^{2+\d, 1+\d / 2}(\mathbb{R} \times[\tau, S])} \leq C_\tau.
\end{align*}
Since $p_{xx}=\phi_{xxx}$, this gives
\begin{align*}
\sup _{\tau \leq t<S}\left\|\phi_{xxx}(\cdot, t)\right\|_{\infty} \leq C_\tau.
\end{align*}
At each step, the previously obtained estimates control the coefficients and the lower order terms in Hölder norm, allowing us to apply Schauder theory to the differentiated equation. Differentiating the equation again, the coefficients and lower order terms involve derivatives already controlled. Repeating the Schauder estimates on successively smaller positive time intervals yields
\begin{align*}
\sup _{\tau \leq t<T_{\max }}\left\|\p_x^m \phi(\cdot, t)\right\|_{\infty} \leq C_{m, \tau},
\end{align*}
for every $m \geq 0$. The constants remain independent of the upper endpoint $S<T_{\max }$.
\end{proof}

Next we show uniqueness. 

\begin{lemma}[Uniqueness from the initial height]\label{lem: uniqueness-height}
Let $T>0$, let $\phi_0\in C_b^1(\R;\R^{1+k})$, and suppose that $K_0 \Subset \SS$. Let $\phi$ and $\psi$ be bounded maps in $C(\R\times[0,T])\cap C^{2,1}(\R\times(0,T])$ with values in $\R^{1+k}$. Suppose each $u\in\{\phi,\psi\}$ satisfies
    \begin{equation*}
        u_t=a(v+u_x)u_{xx},\qquad u(\cdot,0)=\phi_0,\qquad \inf_{\R\times(0,T]}g(v+u_x)>0 .
    \end{equation*}
If for some $\tfrac12<\mu<1$ the solution $\psi$ satisfies
    \begin{equation*}
     M:=\sup_{0<t\le T}t^{1-\mu}\|\psi_{xx}(\cdot,t)\|_\infty<\infty,
    \end{equation*}
then $\phi=\psi$ on $\R\times[0,T]$.
\end{lemma}
\begin{proof}
Let $w=\phi-\psi$, $p=v+\phi_x$ and $q=v+\psi_x$. For each $u\in\{\phi,\psi\}$, uniform spacelikeness gives
\begin{align*}
0<a(v+u_x)
\leq\left(\inf_{\R\times(0,T]}g(v+u_x)\right)^{-1}<\infty,
\end{align*}
so the full height $vx + c + u$ satisfies the hypotheses of \cref{prop: hull preservation} since its perturbation is bounded and continuous up to $t = 0$, and its initial slope is $v + \phi_0'$. Applying \cref{prop: hull preservation} separately to $\phi, \psi$ yields $p, q\in K_0$ for every positive time. Since $K_0 \Subset \SS$ and $a = 1/g$ is smooth and positive on $\SS$, the constants $\lambda = \min_{r \in K_0} a(r) > 0$, $\L = \max_{r\in K_0}a(r) < \infty$ and $L = \max_{r\in K_0} \|Da(r)\| < \infty$ are well defined. The segment joining $p$ to $q$ lies in $K_0$ by convexity, hence the fundamental theorem of calculus along this segment gives 
    \begin{equation*}
        \begin{aligned}
            a(p)-a(q) &=\int_0^1 Da \bigl(q+\t (p-q)\bigr)[p-q]\,d\t,\\
            |a(p)-a(q)|&\leq L\|p-q\|= L\|w_x\|.
        \end{aligned}
    \end{equation*}    
Subtracting the evolution equations for $\phi, \psi$ gives
    \begin{equation*}
        w_t-a(p)w_{xx}=(a(p)-a(q))\psi_{xx}.
    \end{equation*}
We estimate the squared norm of $w$ in order to apply the scalar maximum principle. The fixed positive definite inner product satisfies
    \begin{equation*}
        \p_t\|w\|^2 = 2\llangle w,w_t\rrangle, \qquad \p_{xx}\|w\|^2 = 2\|w_x\|^2 + 2\llangle w,w_{xx}\rrangle.
    \end{equation*}
Consequently, the Cauchy--Schwarz inequality and the preceding bound give
    \begin{equation*}
        \begin{aligned}
            (\p_t-a(p)\p_{xx})\|w\|^2 &=-2a(p)\|w_x\|^2
             +2(a(p)-a(q))\llangle w,\psi_{xx}\rrangle\\
             & \leq - 2 \lambda\|w_x\|^2 + 2L\|w_x\|\,\|w\|\,\|\psi_{xx}\|\\
             & \leq \frac{L^2}{\lambda}\|\psi_{xx}\|^2 \|w\|^2\\
             & \leq At^{2\mu - 2}\|w\|^2,
        \end{aligned}
    \end{equation*}
where the second last inequality follows from Young's inequality, and $A = L^2M^2/\lambda$. 

The assumption $\mu>1/2$ is what makes its upper bound integrable at the initial time:
    \begin{equation*}
        J(t):=\int_0^t A s^{2\mu-2}\,ds =\frac{A}{2\mu-1}t^{2\mu-1}, \qquad J(0)=0.
    \end{equation*}
    Let $U(x,t)=e^{-J(t)}\|w(x,t)\|^2$. For $t>0$,
    \begin{equation*}
        (\p_t-a(p)\p_{xx})U = e^{-J(t)}\left((\p_t-a(p)\p_{xx})\|w\|^2 -J'(t)\|w\|^2\right)\leq 0.
    \end{equation*}
Moreover, $U$ is bounded and continuous on $\R\times[0,T]$, it belongs to $C^{2,1}$ for positive time, and has $U(x,0)=0$ because the initial
heights agree. Since $0<a(p)\leq \L$, we may apply
\cref{lem : maximum principle for slope preservation} with $B=0$, which gives $U\leq 0$ throughout the slab. By definition $U \geq 0$, so $w=0$ and hence $\phi=\psi$.
\end{proof}

We are now ready to prove long time existence and uniqueness of the flow. 

\begin{proposition}[Global existence and uniqueness]
\label{longtimeexistence}
Let $0<\alpha<1$ and suppose $\phi_0\in C_b^{1,\alpha}(\R;\R^{1+k})$. Assume $K_0\Subset\mathcal S$. Then there exists a unique global bounded graphical solution of
\begin{align*}
\phi_t=a(v+\phi_x)\phi_{xx},\qquad \phi(\cdot,0)=\phi_0,
\end{align*}
smooth for positive times and continuous into $C_b^1$ at $t=0$. This solution is unique among bounded classical graphical solutions that are uniformly spacelike on every finite time interval.
\end{proposition}
\begin{proof}
Suppose, for contradiction, that $T_{\max }<\infty$. By \cref{prop: global estimates}, all spatial derivatives of $\phi$ are uniformly bounded on $[T_{\max}/2,T_{\max})$. Differentiating the equation $m$ times gives
\begin{align*}
\p_t \p_x^m \phi=\p_x^m\left(a(p) \phi_{x x}\right) .
\end{align*}
Every term on the right is bounded. Therefore,
\begin{align*}
\sup _{T_{\max } / 2 \leq t<T_{\max }}\left\|\p_t \p_x^m \phi(\cdot, t)\right\|_{\infty}<\infty .
\end{align*}
For $T_{\max}/2\le t_1<t_2<T_{\max }$, from the fundamental theorem of calculus it follows that
\begin{align*}
\left\|\p_x^m \phi\left(\cdot, t_2\right)-\p_x^m \phi\left(\cdot, t_1\right)\right\|_{\infty} \leq C_m\left|t_2-t_1\right| .
\end{align*}
Consequently, $\phi(\cdot,t)$ converges in $C_b^m$ for every integer $m\geq0$ to a limit $\phi_*$ with bounded derivatives of all orders. Since convergence holds in $C_b^2$, it also holds in $C_b^{1+\beta}$ by interpolation. Moreover, $\phi_*$ has bounded derivatives
of all orders and hence belongs to $h_b^{1+\beta}$. Since
$h_b^{1+\beta}$ is closed in $C_b^{1+\beta}$, the convergence also holds in $h_b^{1+\beta}$. The continuation criterion in \cref{lem: local existence} therefore extends the solution beyond $T_{\max}$, contradicting maximality.

Thus $T_{\max}=\infty$. For every finite $T>0$, the weighted estimate in \cref{eq: local-weight} near $t=0$ together with the positive-time derivative bounds of \cref{prop: global estimates} verifies the hypothesis of \cref{lem: uniqueness-height} for the constructed solution, hence uniqueness follows.
\end{proof}

\section{Asymptotic convergence to an affine spacelike line}
\begin{theorem}\label{thm: convergence}
Under the assumptions of \cref{longtimeexistence}, assume in addition that $\phi_0\in C_0(\R;\R^{1+k})$. Let $\phi$ be the unique global solution obtained there. Then, for every integer $m\geq0$,
\begin{align*}
\lim_{t\to\infty}\|\partial_x^m\phi(\cdot,t)\|_\infty=0.
\end{align*}
Consequently, the evolving graphs $F(x,t)=(x,vx+c+\phi(x,t))$ converge smoothly and uniformly,
without rescaling, to the affine spacelike line
\begin{align*}
\ell_{v,c}=\{(x,vx+c):x\in\R\}.
\end{align*}
\end{theorem}
\begin{proof}
We will first prove that $\|\phi(\cdot, t)\|_{\infty} \rightarrow 0$. \cref{prop: hull preservation} gives
	\begin{align*}
p=v+\phi_x \in K_0 \Subset \mathcal{S}.
	\end{align*}
From \cref{phi0infty} we know $\|\phi(\cdot, t)\|_{\infty} \leq\left\|\phi_0\right\|_{\infty}$. However, this last inequality only tells us that $\phi$ remains bounded. It does not imply that $\phi$ converges to zero. We therefore need an additional argument. The key assumption here is $\phi_0 \in C_0$. By definition,
	\begin{align*}
\lim _{|x| \rightarrow \infty} \phi_0(x)=0.
	\end{align*}
Let $\eta>0$. We can choose $R \geq 1$ sufficiently large that
	\begin{align*}
\|\phi_0(x)\| \leq \eta, \quad \text { whenever } |x| \geq R.
	\end{align*}
Thus, the initial perturbation may be large inside $[-R, R]$, but it is small outside this interval. Our aim is to control the part inside $[-R, R]$ by a function that decays to zero as time increases. We will construct a comparison function that controls the potentially large perturbation inside $[-R,R]$ and converges uniformly to zero as $t\to\infty$. Set
	\begin{align*}
M=\max \left\{1,\left\|\phi_0\right\|_{\infty}\right\}, \quad b=\frac{1}{4 \Lambda}, \quad \theta=\frac{\lambda}{2 \Lambda}>0
	\end{align*}
and define
	\begin{align*}
W(x, t)=M e^b\left(\frac{R^2}{t+R^2}\right)^\theta \exp \left(-\frac{b x^2}{t+R^2}\right).
	\end{align*}
We choose $W$ because it satisfies two properties: it controls the potentially large part of the initial perturbation and it converges uniformly to zero as $t \rightarrow \infty$. First, before using the maximum principle, we must verify that our comparison function bounds the initial data. At $t=0$,
	\begin{align*}
W(x, 0)=M e^b \exp \left(-\frac{b x^2}{R^2}\right).
	\end{align*}
For $|x| \leq R$,
	\begin{align*}
W(x, 0) \geq M e^b e^{-b}=M.
	\end{align*}
So, since $\|\phi_0(x)\| \leq M$, we have $W(x, 0) \geq\|\phi_0(x)\|$, inside $[-R, R]$. Outside $[-R, R]$, we already know that $\|\phi_0\| \leq \eta$. Therefore, for every component $\phi_0^i$,
	\begin{align}\label{rangeofphii0}
-\eta-W(x, 0) \leq \phi_0^i(x) \leq \eta+W(x, 0).
	\end{align}
We next show that these inequalities remain true under the evolution. Let $s=t+R^2$. Differentiating gives
	\begin{align*}
W_t=W\left(-\frac{\theta}{s}+\frac{b x^2}{s^2}\right)\quad\text{and}\quad W_{x x}=W\left(-\frac{2 b}{s}+\frac{4 b^2 x^2}{s^2}\right).
	\end{align*}
From \cref{rangeofa}, we have
	\begin{align*}
\frac{\left(\p_t-a(p) \p_{x x}\right) W}{W} & =\frac{-\theta+2 b a(p)}{s}+\frac{\left(b-4 b^2 a(p)\right) x^2}{s^2} \\
& =\frac{a(p)-\lambda}{2 \Lambda s}+\frac{(\Lambda-a(p)) x^2}{4 \Lambda^2 s^2} \\
& \geq 0.
	\end{align*}
Hence, $W$ is a supersolution. Now, every component satisfies
    \begin{align*}
\phi_t^i=a(p) \phi_{x x}^i.
    \end{align*}
We want to show that \cref{rangeofphii0} holds for all times. Since $W$ is a supersolution, the functions
	\begin{align}\label{2functions}
\phi^i-\eta-W \quad \text { and } \quad-\phi^i-\eta-W
	\end{align}
are scalar subsolutions with nonpositive initial values. Indeed, for $f=\phi^i-\eta-W$, notice that
	\begin{align*}
f_t-a(p) f_{x x} & =\left(\phi_t^i-a(p) \phi_{x x}^i\right)-\left(W_t-a(p) W_{x x}\right) \\
& =0-\left(W_t-a(p) W_{x x}\right) \\
& \leq 0 .
	\end{align*}
Therefore, $f$ is a subsolution. Moreover, $f(x, 0) \leq 0$. Applying \cref{lem : maximum principle for slope preservation} to each of the two functions in \cref{2functions} gives the desired comparison inequalities
	\begin{align*}
\left|\phi^i(x, t)\right| \leq \eta+W(x, t).
	\end{align*}
Taking the supremum over $x$ and observing that
	\begin{align*}
\sup _{x \in \mathbb{R}} W(x, t)=M e^b\left(\frac{R^2}{t+R^2}\right)^\theta,
	\end{align*}
we obtain
	\begin{align*}
\left\|\phi^i(\cdot, t)\right\|_{\infty} \leq \eta+M e^b\left(\frac{R^2}{t+R^2}\right)^\theta.
	\end{align*}
Since $\theta>0$, for the fixed $R$ chosen above,
	\begin{align*}
\left(\frac{R^2}{t+R^2}\right)^\theta \longrightarrow 0, \quad \text { as } t \rightarrow \infty.
	\end{align*}
Therefore,
	\begin{align*}
\limsup _{t \rightarrow \infty}\left\|\phi^i(\cdot, t)\right\|_{\infty} \leq \eta.
	\end{align*}
But $\eta>0$ was arbitrary. Letting $\eta \rightarrow 0$, we conclude that
	\begin{align*}
\left\|\phi^i(\cdot, t)\right\|_{\infty} \longrightarrow 0.
	\end{align*}
There are only finitely many components, so 
	\begin{align}\label{phito0}
\|\phi(\cdot, t)\|_{\infty} \longrightarrow 0.
	\end{align}
We now prove convergence of all spatial derivatives. By the uniform higher derivative estimates in \cref{derest}, for every $m\geq0$ and $\tau>0$, we have
	\begin{align*}
\sup _{t \geq \tau}\left\|\p_x^m \phi(\cdot, t)\right\|_{\infty} \leq C_{m, \tau}.
	\end{align*}
We want to strengthen them to convergence to zero. To do this, we use the elementary interpolation inequality
	\begin{align}\label{interpolationeq}
\left\|f^{\prime}\right\|_{\infty}^2 \leq 4\|f\|_{\infty}\left\|f^{\prime \prime}\right\|_{\infty}
	\end{align}
for a bounded twice-differentiable function with bounded second derivative. Let's recall why. By Taylor's formula,
	\begin{align*}
f(x+h)=f(x)+h f^{\prime}(x)+\int_0^h(h-s) f^{\prime \prime}(x+s)ds.
	\end{align*}
Rearranging and taking norms,
	\begin{align*}
\|f^{\prime}(x)\| \leq \frac{2\|f\|_{\infty}}{h}+\frac{h}{2}\|f^{\prime \prime}\|_{\infty}.
	\end{align*}
If $\|f''\|_\infty=0$, then $f$ is a bounded affine function and hence constant, so the inequality is immediate. Choosing $h=2 \sqrt{\frac{\|f\|_{\infty}}{\left\|f^{\prime \prime}\right\|_{\infty}}}$ to minimise the right-hand side gives \cref{interpolationeq}. We first apply it to $f=\phi(\cdot, t)$ to get
	\begin{align*}
\left\|\phi_x(\cdot, t)\right\|_{\infty}^2 \leq 4\|\phi(\cdot, t)\|_{\infty}\left\|\phi_{xx}(\cdot, t)\right\|_{\infty}.
	\end{align*}
From \cref{phito0} and \cref{boundedphixx}, we have
	\begin{align}\label{phixto0}
\left\|\phi_x(\cdot, t)\right\|_{\infty} \longrightarrow 0.
	\end{align}
Next, we apply \cref{interpolationeq} to $f=\phi_x(\cdot, t)$ to get
	\begin{align*}
\left\|\phi_{xx}(\cdot, t)\right\|_{\infty}^2 \leq 4\left\|\phi_x(\cdot, t)\right\|_{\infty}\left\|\phi_{xxx}(\cdot, t)\right\|_{\infty} .
	\end{align*}
From \cref{derest} and \cref{phixto0}, we get
	\begin{align*}
\left\|\phi_{xx}(\cdot, t)\right\|_{\infty} \longrightarrow 0.
	\end{align*}
We repeat this argument. More precisely, suppose that
	\begin{align*}
\left\|\p_x^{m-1} \phi(\cdot, t)\right\|_{\infty} \longrightarrow 0.
	\end{align*}
The interpolation inequality gives
\begin{align*}
\left\|\p_x^m \phi(\cdot,t)\right\|_\infty^2&\leq 4\left\|\p_x^{m-1}\phi(\cdot,t)\right\|_\infty\cdot\left\|\p_x^{m+1}\phi(\cdot,t)\right\|_\infty.
\end{align*}
The first factor tends to zero by the induction hypothesis, while the second is uniformly bounded by \cref{derest}. Therefore,
	\begin{align*}
\left\|\p_x^m \phi(\cdot, t)\right\|_{\infty} \longrightarrow 0.
	\end{align*}
By induction, we conclude
	\begin{align}\label{finalderest}
\lim _{t \rightarrow \infty}\left\|\p_x^m \phi(\cdot, t)\right\|_{\infty}=0, \quad \text { for every } m \geq 0.
	\end{align}
Finally, since $F(x, t)=(x, v x+c+\phi(x, t))$, we have
	\begin{align*}
F(x, t)-(x, v x+c)=(0, \phi(x, t))
	\end{align*}
and using \cref{finalderest} shows that the graphs converge uniformly, together with all their spatial derivatives, to the affine spacelike line $\ell_{v,c}$. This completes the proof. 
\end{proof}

We conclude this section by giving the proof of Theorem~\ref{thm: main theorem}.

\begin{proof}[Proof of Theorem~\ref{thm: main theorem}]
Global existence and uniqueness follow from
Proposition~\ref{longtimeexistence}. By Proposition~\ref{prop: hull preservation},
\begin{align*}
K_t\subseteq K_s\subseteq K_0\Subset\mathcal S,\qquad 0\leq s\leq t.
\end{align*}
Finally, when $\phi_0\in C_0$, Theorem~\ref{thm: convergence} gives smooth and uniform convergence to $\ell_{v,c}$, completing the proof.
\end{proof}

\section{Ancient graphical rigidity}

In this section, we give the proof of Theorem~\ref{thm:ancient-rigidity}. The argument uses the uniform parabolicity of the graphical equation and the interior H\"older estimates for scalar divergence-form parabolic equations.

\begin{proof}[Proof of Theorem~\ref{thm:ancient-rigidity}]
By the graphical formulation of curve shortening flow in \cref{graphical MCF equation}, the graph function satisfies
\begin{align*}
u_t=a(u_x)u_{xx},\qquad a(p)=\frac{1}{g(p)}.
\end{align*}
Setting $p=u_x$ and differentiating gives
\begin{align}\label{eq:ancient-slope}
p_t=\partial_x\big(a(p)p_x\big).
\end{align}
Indeed, the mean curvature vector is the normal projection of $(0,u_{xx}/g(u_x))$. Since the first
coordinate of $F_x=(1,u_x)$ is nonzero, the graphical velocity $(0,u_t)$ has the same normal projection precisely when $u_t=u_{xx}/g(u_x)$.

We will show that $u_x$ is constant in both
space and time. Set $b(x,t):=a(p(x,t))=\frac{1}{g(p(x,t))}$. Although $p$ is vector-valued, each component $p^j$ satisfies the scalar divergence-form parabolic equation
\begin{align}\label{eq:ancient-scalar}
\partial_t p^j
=\partial_x\bigl(b(x,t)\partial_xp^j\bigr).
\end{align}
By assumption, there are constants $M<\infty$ and $\delta>0$, such that
\begin{align*}
\|p(x,t)\|\leq M,\qquad g(p(x,t))\geq\delta
\end{align*}
throughout $\R\times(-\infty,T)$. Since
\begin{align*}
g(p)=1+(p^+)^2-\sum_{j=1}^k(p^j)^2\leq 1+\|p\|^2,
\end{align*}
we obtain
\begin{align}\label{eq:ancient-ellipticity}
0<\lambda:=\frac{1}{1+M^2}\leq b(x,t)\leq\frac{1}{\delta}=:\Lambda<\infty.
\end{align}
Thus \cref{eq:ancient-scalar} is uniformly parabolic, with ellipticity constants independent of space and time.

We now apply Nash's interior H\"older estimate \cite[Sections 4,6]{Nash} for scalar divergence-form parabolic equations. Recall that if a bounded function $w$ solves
\begin{align*}
w_s=\partial_y\bigl(B(y,s)w_y\bigr)
\end{align*}
on the parabolic cylinder
\begin{align*}
Q_1=(-1,1)\times(-1,0],
\end{align*}
where $\lambda\leq B\leq\Lambda$, then there exist constants $C<\infty$ and $\gamma\in(0,1)$, depending only on $\lambda$ and $\Lambda$, such that
\begin{align}\label{eq:ancient-nash}
|w(y,s)-w(y',s')|\leq C\|w\|_{L^\infty(Q_1)}\left(|y-y'|+|s-s'|^{1/2}\right)^\gamma
\end{align}
for all $(y,s),(y',s')\in Q_{1/2}$, where
\begin{align*}
Q_{1/2}=(-1/2,1/2)\times(-1/4,0].
\end{align*}
Importantly, the constants depend only on the ellipticity bounds, not on derivatives of $B$.

Fix two arbitrary points
\begin{align*}
(x_1,t_1),\ (x_2,t_2)
\in\R\times(-\infty,T).
\end{align*}
Choose $x_*\in\R$ and $\max\{t_1,t_2\}<t_*<T$. For each $R>0$, define the rescaled slope components
\begin{align*}
w_R^j(y,s)=p^j(x_*+Ry,t_*+R^2s),
\qquad (y,s)\in Q_1.
\end{align*}
This rescaling is well-defined for every $R>0$ because the original solution is defined for all times less than $T$.

By \cref{eq:ancient-scalar}, each $w_R^j$ satisfies
\begin{align*}
\partial_s w_R^j=\partial_y\left(B_R(y,s)\partial_yw_R^j\right),
\end{align*}
where
\begin{align*}
B_R(y,s)=b(x_*+Ry,t_*+R^2s).
\end{align*}
The bounds in \cref{eq:ancient-ellipticity} are unchanged by this rescaling
\begin{align*}
\lambda\leq B_R(y,s)\leq\Lambda,\qquad\|w_R^j\|_{L^\infty(Q_1)}\leq M.
\end{align*}
For sufficiently large $R$, the two rescaled points
\begin{align*}
(y_i,s_i)=\left(\frac{x_i-x_*}{R}, \frac{t_i-t_*}{R^2}\right),\qquad i=1,2,
\end{align*}
belong to $Q_{1/2}$. Applying \cref{eq:ancient-nash} gives
\begin{align*}
|p^j(x_1,t_1)-p^j(x_2,t_2)|\leq CM\left(\frac{|x_1-x_2|}{R}+\frac{|t_1-t_2|^{1/2}}{R}\right)^\gamma.
\end{align*}
The constants $C,M,\gamma$ are independent of $R$. Sending $R\to\infty$, we conclude
\begin{align*}
p^j(x_1,t_1)=p^j(x_2,t_2).
\end{align*}
Since both points and the component $j$ were arbitrary, every component of $p$ is constant on $\R\times(-\infty,T)$. Therefore
\begin{align*}
u_x(x,t)=p_*,
\end{align*}
for some fixed vector $p_*\in\R^{1+k}$. Integrating with respect to $x$, we obtain
\begin{align*}
u(x,t)=p_*x+c(t),
\end{align*}
for some smooth function $c:(-\infty,T)\to\R^{1+k}$. In particular, $u_{xx}=0$. Substituting into
\cref{graphical MCF equation} gives
\begin{align*}
c'(t)=u_t=a(p_*)u_{xx}=0.
\end{align*}
Hence $c(t)=c_*$ is constant and
\begin{align*}
F(x,t)=(x,p_*x+c_*)
\end{align*}
for all $x\in\R$ and $t<T$.

Finally, \cref{eq:ancient-bounds} implies
$g(p_*)\geq\delta>0$, so the affine
line is spacelike. This completes the proof.
\end{proof}

\bibliography{bib}

@misc{andrews2026highcodimensionmeancurvature,
      title={High codimension mean curvature flow of spacelike-convex submanifolds with one spacelike codimension}, 
      author={Ben Andrews and Qiyu Zhou},
      year={2026},
      howpublished={https://arxiv.org/abs/2602.19537}, 
}

@article{AndrewsBryan2011,
  author  = {Andrews, Ben and Bryan, Paul},
  title   = {Curvature bound for curve shortening flow via distance comparison and a direct proof of Grayson’s theorem},
  journal = {Journal für die reine und angewandte Mathematik},
  volume  = {653},
  pages   = {179–187},
  year    = {2011},
  doi     = {10.1515/crelle.2011.026}
}

@book{BoydVandenberghe,
  author    = {Stephen Boyd and Lieven Vandenberghe},
  title     = {Convex Optimization},
  publisher = {Cambridge University Press},
  year      = {2004},
  url       = {https://web.stanford.edu/~boyd/cvxbook/bv_cvxbook.pdf}
}

@article{ClementSimonett,
  author  = {Cl{\'e}ment, Philippe and Simonett, Gieri},
  title   = {Maximal regularity in continuous interpolation spaces and quasilinear parabolic equations},
  journal = {Journal of Evolution Equations},
  year    = {2001},
  volume  = {1},
  number  = {1},
  pages   = {39--67},
  doi     = {10.1007/PL00001364}
}

@article{ShaoSimonett,
  author  = {Shao, Yuanzhen and Simonett, Gieri},
  title   = {Continuous maximal regularity on uniformly regular {Riemannian} manifolds},
  journal = {Journal of Evolution Equations},
  year    = {2014},
  volume  = {14},
  number  = {1},
  pages   = {211--248},
  doi     = {10.1007/s00028-014-0218-6}
}

@book{LSU,
  author    = {Ladyzhenskaya, O. A. and Solonnikov, V. A. and Ural'tseva, N. N.},
  title     = {Linear and Quasi-linear Equations of Parabolic Type},
  series    = {Translations of Mathematical Monographs},
  volume    = {23},
  publisher = {American Mathematical Society},
  address   = {Providence, RI},
  year      = {1968},
  note      = {Translated from the Russian by S. Smith},
  doi       = {10.1090/mmono/023}
}

@article{Nash,
  author    = {Nash, John},
  title     = {Continuity of Solutions of Parabolic and Elliptic Equations},
  journal   = {American Journal of Mathematics},
  volume    = {80},
  number    = {4},
  pages     = {931--954},
  year      = {1958},
  doi       = {10.2307/2372841}
}

@article{LambertLotay,
  author    = {Lambert, Ben and Lotay, Jason D.},
  title     = {Spacelike Mean Curvature Flow},
  journal   = {Journal of Geometric Analysis},
  volume    = {31},
  pages     = {1291--1359},
  year      = {2021},
  doi       = {10.1007/s12220-019-00266-4}
}

@article{GageHamilton,
  author  = {Gage, Michael E. and Hamilton, Richard S.},
  title   = {The Heat Equation Shrinking Convex Plane Curves},
  journal = {Journal of Differential Geometry},
  volume  = {23},
  number  = {1},
  pages   = {69--96},
  year    = {1986},
  doi     = {10.4310/jdg/1214439902},
  url     = {https://doi.org/10.4310/jdg/1214439902}
}

@article{Grayson,
  author  = {Grayson, Matthew A.},
  title   = {The Heat Equation Shrinks Embedded Plane Curves to Round Points},
  journal = {Journal of Differential Geometry},
  volume  = {26},
  number  = {2},
  pages   = {285--314},
  year    = {1987},
  doi     = {10.4310/jdg/1214441371},
  url     = {https://doi.org/10.4310/jdg/1214441371}
}

@article{EckerHuisken,
  author  = {Ecker, Klaus and Huisken, Gerhard},
  title   = {Mean Curvature Evolution of Entire Graphs},
  journal = {Annals of Mathematics},
  volume  = {130},
  number  = {3},
  pages   = {453--471},
  year    = {1989},
  doi     = {10.2307/1971452},
  url     = {https://doi.org/10.2307/1971452}
}

@article{LiSalavessa,
  author  = {Li, Guanghan and Salavessa, Isabel M. C.},
  title   = {Mean Curvature Flow of Spacelike Graphs},
  journal = {Mathematische Zeitschrift},
  volume  = {269},
  number  = {3--4},
  pages   = {697--719},
  year    = {2011},
  doi     = {10.1007/s00209-010-0768-4},
  url     = {https://doi.org/10.1007/s00209-010-0768-4}
}

@article{MiuraRupp,
  author  = {Miura, Tatsuya and Rupp, Fabian},
  title   = {A New Energy Method for Shortening and Straightening Complete Curves},
  journal = {Journal f{\"u}r die reine und angewandte Mathematik},
  volume  = {2026},
  number  = {830},
  pages   = {1--50},
  year    = {2026},
  doi     = {10.1515/crelle-2025-0077},
  url     = {https://doi.org/10.1515/crelle-2025-0077}
}

@misc{HuMa,
  author        = {Hu, Baichuan and Ma, Xiang},
  title         = {Weighted Normalized Curve Shortening Flow with Applications in Pseudo-{Euclidean} Spaces},
  year          = {2026},
  eprint        = {2609.24027},
  archivePrefix = {arXiv},
  primaryClass  = {math.DG},
  howpublished  = {https://arxiv.org/abs/2609.24027}
}

@misc{NguyenVogiatziFBCSF,
  author        = {Nguyen, Huy The and Vogiatzi, Artemis A.},
  title         = {High Codimension Curve Shortening Flow with Free Boundary},
  year          = {2026},
  howpublished =  {arXiv:2602.20865}
}

@article{AltschulerGrayson1992,
  author  = {Altschuler, Steven J. and Grayson, Matthew A.},
  title   = {Shortening Space Curves and Flow Through Singularities},
  journal = {Journal of Differential Geometry},
  volume  = {35},
  number  = {2},
  pages   = {283--298},
  year    = {1992},
  doi     = {10.4310/jdg/1214448076}
}
  \bibliographystyle{abbrv}

\end{document}